\documentclass{article}

\usepackage{amsmath,amssymb,graphicx} 
\usepackage[margin=1in]{geometry}
\usepackage{enumerate}
\begin{document}

\nocite{*} 

\title{Irreversible homotopy and a notion of irreversible Lusternik–Schnirelmann category}

\author{Khashayar Rahimi }

\date{October 21, 2020} 

\maketitle

\begin{abstract}
This work was intended as an attempt to investigate a model of irreversible process and natural phenomena. For this, we introduce the notion of irreversible path (that for brevity we write ir-path), ir-homotopy, ir-contractible space, and Lusternik-Schnirelmann ir-category by equipping the $I=[0,1]$ with left order topology. We will restrict the irreversibility of definitions to $T_0$ Spaces, such that for $T_1$ spaces, the ir-paths are constant. After providing some theorems and properties of these notions, eventually, we prove that Lusternik-Schnirelmann ir-category is an invariant of ir-homotopy equivalence.
\end{abstract}

\section{Introduction}
In this note, we wish to introduce another notion of path, homotopy, contractibility and homotopy equivalence that are not necessarily reversible, than the one of [1],[2].\\

\textbf{Definition 1.} The collection  $\mathcal{B} = \{ (- \infty,b) ; b \in \mathbb{R}  \}$ of subsets of $\mathbb{R}$ is a basis for a topology on $\mathbb{R}$. If we denote this topology by $\tau$, then  $\tau = \mathcal{B} \cup \{\emptyset , \mathbb{R}\}$. The topological space $(\mathbb{R},\tau)$ denoted by $\stackrel{\rm ir}{\rm \mathbb{R}}$. The subspace topology of subset $I = [0,1]$ of $\stackrel{\rm ir}{\rm \mathbb{R}}$, is the collection $\{[0,b) ; 0 \le b \le 1    \} \cup \{I\}$, and denoted by $\stackrel{\rm ir}{\rm I}$.\
The spaces  $\mathbb{R}$ and  $I$ with standard topology, respectively denoted by $\mathbb{R}$ and  $I$.\\

$\tau$ is called left order topology on $\mathbb{R}$. [3]

\vspace{1cm}
\textbf{Proposition 1.} Let $I = [0,1]$ then $d_{I}^{ir} : I \times I \to \overline{\mathbb{R}^{+}}$ where $d_{I}^{ir} = max \{y-x,0\}$ is a quasi-metric [4] and its open ball topology $\mathcal{O}^{d_{I}^{ir}}$ is equal to left order topology on $I$.

\vspace{1cm}
\textit{Proof.} It's easy to check the three conditions of quasi-metric spaces for ${d_{I}}^{ir}$. Also, the open ball $B_{x , < \epsilon }^{d_{I}^{ir}}$ with center $x$ and radius $\epsilon \in \mathbb{R}^{+}$ is of the following form\\
\begin{center}
$B_{x , < \epsilon}^{d_{I}^{ir}} = \bigg \{ y \in I  \ \bigg | \ max \{y-x,0\} < \epsilon \bigg \} = [0, x+ \epsilon)$
\end{center}

Thus, the open balls of $I$ are equal to open sets of left order topology on $I$, and consequently $\mathcal{O}^{{d_{I}}^{ir}}$ is equal to $\stackrel{\rm ir}{\rm I}$.

\begin{flushright}
$\square$
\end{flushright}

\textbf{Proposition 2.} Every subspace of $\stackrel{\rm ir}{\rm \mathbb{R}}$ is hyperconnected. In particular, $\stackrel{\rm ir}{\rm \mathbb{R}}$ and $\stackrel{\rm ir}{\rm I}$ are connected.

\vspace{1cm}

\textit{Proof.} The intersection of every two nonempty and open subset in $\stackrel{\rm ir}{\rm \mathbb{R}}$ is a nonempty subset.

\begin{flushright}
$\square$
\end{flushright}

\textbf{Theorem 1.} A subset $A$ of $\stackrel{\rm ir}{\rm \mathbb{R}}$ is compact if and only if it has a biggest element. 

\vspace{1cm}

\textit{Proof.} Suppose that $m$ be the biggest element of $A$ and  $\mathcal{U}$ be a cover of open subsets of  $\stackrel{\rm ir}{\rm \mathbb{R}}$. An open subset $U$ of $\mathcal{U}$ exists such that  $m \in U$. Thus $(- \infty , m] \subseteq U$ and so $A$ is compact.\\
Conversely, assume that $A$ is a subset of $\stackrel{\rm ir}{\rm \mathbb{R}}$ that dosen't have the biggest element. Thus the collection $\mathcal{V} = \{(- \infty ,a) ; a \in A  \}$ is an open cover for $A$ that doesn't have any finite subcover, then $A$ is not compact.

\begin{flushright}
$\square$
\end{flushright}

\textbf{Corollary 1.} $\stackrel{\rm ir}{\rm \mathbb{R}}$ is not and $\stackrel{\rm ir}{\rm I}$ is compact.

\section{ir-path}

\textbf{Definition 2.} Let $X$ be a topological space. A function  $\gamma : \stackrel{\rm ir}{\rm I} \to X$ is called an ir-path in $X$, if it is continuous on $\stackrel{\rm ir}{\rm I}$. $\gamma(0)$ is the initial point and  $\gamma(1)$ is the terminal point of the  ir-path $\gamma$. 

\vspace{1cm}

\textbf{Theorem 2.} Let $X$ be a topological space and $x,y$ two points in $X$. Then an ir-path from $x$ to $y$ exists, if and only if $y \in \overline{\{x\}}$.
\vspace{1cm}
\\
\textit{Proof.} Suppose that $y \in \overline{\{x\}}$. Then the function $\gamma : \stackrel{\rm ir}{\rm I} \to X$ with the following formula is an ir-path from $x$ to $y$. 
\begin{center}
$\gamma (s) =
\bigg \{
	\begin{array}{ll}
		x  &  \hspace{1cm}  0 \le t < 1 \\
		y &  \hspace{1cm} t = 1 \\
	\end{array}$
\end{center}
For proving, assume that $V$ is an open subset of $X$. If $x,y \notin V$ then ${\gamma}^{-1} (V) = \emptyset$ that is an open subset of $\stackrel{\rm ir}{\rm I}$. Now if $y \in V$ then $x \in V$, thus ${\gamma}^{-1} (V) = \stackrel{\rm ir}{\rm I}$ that is an open subset of $\stackrel{\rm ir}{\rm I}$. Also, if $y \notin V$ and $x \in V$, then ${\gamma}^{-1} (V) = [0,1)$ is an open subset of $\stackrel{\rm ir}{\rm I}$. Thus, $\gamma : \stackrel{\rm ir}{\rm I} \to X$ is an ir-path from $x$ to $y$.\\
Conversely, suppose that $\gamma : \stackrel{\rm ir}{\rm I} \to X$ is an ir-path from $x$ to $y$. Assume that $s$ is an arbitrary element of $\stackrel{\rm ir}{\rm I}$ and $z = \gamma(s)$. If $V$ is an open subset of $X$ containing $z$,  ${\gamma}^{-1} (V)$  is an open subset of $ \stackrel{\rm ir}{\rm I}$ containing $s$ and also $\{0\}$. Thus $x = \gamma (0) \in V$ then $z \in \overline{\{x\}}$  and in particular,  $y \in \overline{\{x\}}$. 

\begin{flushright}
$\square$
\end{flushright}
\textbf{Theorem 3.} Let $X$ be a topological space and $x \in X$. If $\gamma : \stackrel{\rm ir}{\rm I} \to X$ is an ir-path with $x$ as its initial point, then $\gamma(\stackrel{\rm ir}{\rm I}) \subseteq \overline{\{x\}}$.

\vspace{1cm}

\textit{Proof.} The second part of the proof of Theorem 2.

\begin{flushright}
$\square$
\end{flushright}
\textbf{Theorem 4.} If $X$ be a $T_1$ topological space, then each ir-path in $X$ is constant.

\vspace{1cm}

\textit{Proof.} Suppose that $\gamma : \stackrel{\rm ir}{\rm I} \to X$ is an ir-path in $X$. For every $s \in \stackrel{\rm ir}{\rm I}$, the subset $\{\gamma (s)\}$ of $X$ is closed. Thus from continuity of $\gamma$, the subset ${\gamma}^{-1} \big ( \{\gamma (s)\} \big )$ of $\stackrel{\rm ir}{\rm I}$ is also closed. On the other hand, every closed and nonempty subset of $\stackrel{\rm ir}{\rm I}$ containing the point $\{1\}$, hence ${\gamma}^{-1} \big ( \{\gamma (s)\} \big )$ contains the point $\{1\}$. Thus $\gamma (s) = \gamma (1)$.

\begin{flushright}
$\square$
\end{flushright}

\section{ir-homotopy and ir-contractibility}

\textbf{Definition 3.} Let $X$ and $Y$ be two topological spaces, and $f: X \to Y$ and $g: X \to Y$ be continuous functions on $X$. The function $F: X \times \stackrel{\rm ir}{\rm I} \to Y$ is called an ir-homotopy from $f$ to $g$, if for every $x \in X$ we have $F(x,0) = f(x)$ and $F(x,1) = g(x)$. If there exists an ir-homotopy from $f$ to $g$, we write  $f  \stackrel{\rm ir}{\rm \cong} g$.

\vspace{1cm}

\textbf{Definition 4.} Suppose that $x_0 \in X$. The topological space $X$ is called ir-contractible in $x_0$, if there exists an ir-homotopy from identity function $1_X : X \to X$ to constant function $x_0: X \to X$.

\vspace{1cm}

\textbf{Theorem 5.} Let $X$ and $Y$ be topological spaces and $f: X \to Y$ and $g: X \to Y$ be continuous functions on $X$. If $Y$ be a $T_1$ space and $f  \stackrel{\rm ir}{\rm \cong} g$, then $f=g$.

\vspace{1cm}

\textit{Proof.} Suppose that $F: X \times \stackrel{\rm ir}{\rm I} \to Y$ is an ir-homotopy from $f$ to $g$. For every $x \in X$, we can define a function $F_x: \stackrel{\rm ir}{\rm I} \to Y$ with the formula $F_x(t) = F(x,t)$, such that $F_x$ is an ir-path. Using Theorem 4 , $F_x$ is an constant function. Thus for every $x \in X$ we have: $f(x) = F_x (0) = F_x(1) = g(x)$.

\begin{flushright}
$\square$
\end{flushright}

\vspace{1cm}

\textbf{Notations.} Let $X$ be a topological space. We set $\stackrel{\rm ir}{\rm Co}(X) = \bigcap_{x \in X} \overline{\{x\}}$. In other words, $\stackrel{\rm ir}{\rm Co (X)}$ is the set of all points of $X$ that have just one open neighborhood, and that is $X$. 

\vspace{1cm}

\textbf{Example.}  $\stackrel{\rm ir}{\rm Co} (\stackrel{\rm ir}{\rm \mathbb{R}}) = \stackrel{\rm ir}{\rm Co} (\mathbb{R}) = \stackrel{\rm ir}{\rm Co} (I) = \emptyset$.

\vspace{1cm}

\textbf{Example.} If we consider the Sierpiński space the set $\{0,1\}$ and its topology  $\big \{\emptyset,S ,\{0\}   \big\}$, then\\
 $\stackrel{\rm ir}{\rm Co} (S) = \{1\}$.

\vspace{1cm}

\textbf{Theorem 6.} Let $X$ be a topological space and $x_0 \in X$. The space $X$ is ir-contractible in $x_0$ if and only if $x_0 \in \stackrel{\rm ir}{\rm Co (X)}$.

\vspace{1cm}

\textit{Proof.} Suppose that $x_0 \in \stackrel{\rm ir}{\rm Co (X)}$. We will show that the function $F: X \times \stackrel{\rm ir}{\rm I} \to X$ with the following formula is an ir-homotopy from identity function $1_X : X \to X$ to constant function $x_0: X \to X$:

\begin{center}
$ F(x,t) =
\bigg \{
	\begin{array}{ll}
		x  &  \hspace{1cm}  0 \le t < 1 \\
		x_0 &  \hspace{1cm} t = 1 \\
	\end{array}$
\end{center}
It is clear that for all $x \in X$,  $F(x,0) = x$ and $F(x,1) = x_0$. For proving the continuity of $F$, assume that $V$ is an arbitrary open subset of $X$. If $x_0 \in V$, since $x_0 \in \stackrel{\rm ir}{\rm Co (X)}$ we conclude that $V = X$. Hence,  $F^{-1}(V) = X \times \stackrel{\rm ir}{\rm I}$ is an open subset of $X \times \stackrel{\rm ir}{\rm I}$. Besides, if $x_0 \notin V$, by the formula of $F$, we have $F(V) = V \times [0,1)$. Thus in this case as well as pervious case, $F^{-1}(V)$ is an open subset of $X \times \stackrel{\rm ir}{\rm I}$.\\

Conversely, If $X$ is ir-contractible in $x_0$, There is an ir-homotopy from identity function $1_X : X \to X$ to constant function $x_0: X \to X$. For all $x \in X$, we can define the function $F_x: \stackrel{\rm ir}{\rm I} \to X$ by $F_x(t) = F(x,t)$, that is continuous on $\stackrel{\rm ir}{\rm I}$ and obviousely on the point $\{1\}$. Thus for all arbitrary open neighborhoods $V$ of $x_0$,  ${F_x}^{-1} (V)$ is an open neighborhood of the point $1$ in $\stackrel{\rm ir}{\rm I}$ and therefore it is equal to $\stackrel{\rm ir}{\rm I}$. Now  $0 \in {F_x}^{-1} (V)$ implies that $x = F_x (0) \in V$. But we considered $x$ arbitrarily, hence $V =X$. 
\begin{flushright}
$\square$
\end{flushright}

\vspace{1cm}

\textbf{Corollary 2.}  $\mathbb{R}$,  $I$ and $\stackrel{\rm ir}{\rm \mathbb{R}}$ are not ir-contractible, but the Sierpiński space is ir-contractible.

\vspace{1cm}
\textbf{Corollary 3.} The space $\stackrel{\rm ir}{\rm I}$ is ir-contractible and precisely,  $\stackrel{\rm ir}{\rm Co} (\stackrel{\rm ir}{\rm I}) = \{1\}$. Also, we can consider the following ir-contractiblity function;\\

\begin{center}
$G:\stackrel{\rm ir}{\rm I} \times \stackrel{\rm ir}{\rm I} \to \stackrel{\rm ir}{\rm I} $\\
$G(s,t) = (1-t)s + t$
\end{center} 

\vspace{1cm}

\textbf{Corollary 4.} All the $T_1$ and ir-contractible spaces are singletons.

\vspace{1cm}

\textbf{Corollary 5.} All the ir-contractible spaces are compact.

\vspace{1cm}

\textit{Proof.} Suppose that $X$ is a topological space and it is ir-contractible in $x_0$. If $\mathcal{U}$ be an arbitrary open cover for $X$, there is an open set of $\mathcal{U}$ like $U$, that  $x_0 \in U$. Now by using the Theorem 6, we conclude that $U = X$. Thus, $\mathcal{U}$ has a singleton subcover for $X$.

\begin{flushright}
$\square$
\end{flushright}

\vspace{1cm}

\textbf{Notation.} Let $n$ be a natural number. We denote  $ \underbrace{\mathbb{R} \times \mathbb{R} \times \ldots \times \mathbb{R}}_n$ by $\mathbb{R}^n$, and $\underbrace{\stackrel{\rm ir}{\rm \mathbb{R}} \times \stackrel{\rm ir}{\rm \mathbb{R}} \times \ldots \times \stackrel{\rm ir}{\rm \mathbb{R}}}_n$ by  $\stackrel{\rm ir}{\rm \mathbb{R}^n}$ with their product topology.\\
Besides, we consider $\mathbb{S}^{n-1}$ as the subspace $\{(x_1, \ldots , x_n) \in \mathbb{R}^n : \displaystyle\sum_{i=1}^{n} {x_i}^{2} = 1\} $  of $\mathbb{R}^n$, and  $\stackrel{\rm ir}{\rm \mathbb{S}^{n-1}}$ as the subspace  $\{(x_1, \ldots , x_n) \in \stackrel{\rm ir}{\rm \mathbb{R}^n} :   \displaystyle\sum_{i=1}^{n} {x_i}^{2} = 1\} $ of $\stackrel{\rm ir}{\rm \mathbb{R}^n}$.  

\textbf{Remark.} The topology of $\mathbb{R}^n$ is finer than the topology of $\stackrel{\rm ir}{\rm \mathbb{R}^n}$, thus the topology of  $\mathbb{S}^{n-1}$ is finer than the topology of  $\stackrel{\rm ir}{\rm \mathbb{S}^{n-1}}$.\\
Also, we know that  $\mathbb{S}^{n-1}$ is compact, thus  $\stackrel{\rm ir}{\rm \mathbb{S}^{n-1}}$ is compact too.

\vspace{1cm}

\textbf{Theorem 7.} Suppose that  $(X_\alpha)_{\alpha \in I}$ is a collection of topological spaces. Consider  $(X_\alpha)_{\alpha \in I}$ with product topology. $\prod_{\alpha \in I} X_{\alpha}$ is ir-contractible in $(x_\alpha)_{\alpha \in I}$ if and only if for each  $\alpha$ in $I$, $X_{\alpha}$ is ir-contractible in $x_{\alpha}$.

\vspace{1cm}
\textit{Proof.} Assume that each $X_{\alpha}$ is ir-contractible in $x_{\alpha}$. For each neighborhood $V$ of  $(x_\alpha)_{\alpha \in I}$, there exist open neighborhood $U_{\alpha}$ of $x_{\alpha}$ such that $\prod_{\alpha \in I} U_{\alpha} \subseteq V$. Since each $X_{\alpha}$ is ir-contractible in $x_{\alpha}$, for each $\alpha$ we have  $U_{\alpha} = X_{\alpha}$. Therefore $V = \prod_{\alpha \in I} X_{\alpha}$. Then $\prod_{\alpha \in I} X_{\alpha}$ is ir-contractible in $(x_\alpha)_{\alpha \in I}$.\\

Conversely, suppose that $\prod_{\alpha \in I} X_{\alpha}$ is ir-contractible in $(x_\alpha)_{\alpha \in I}$. Consider an arbitrary $\beta$ of $I$. Assume that $U_{\beta}$ is an arbitrary open neiborhood of  $x_{\beta}$. For each $\alpha$ in $I \setminus \{\beta \}$ we set $U_{\alpha} = X_{\alpha}$. So in this case,  $\prod_{\alpha \in I} U_{\alpha}$ become an open neighborhood of $(x_\alpha)_{\alpha \in I}$. Thus by the assumption we have $\prod_{\alpha \in I} U_{\alpha} = \prod_{\alpha \in I} X_{\alpha}$ and consequently  $U_{\beta} = X_{\beta}$. Therefore  $X_{\beta}$ is ir-contractible in  $x_{\beta}$.   

\begin{flushright}
$\square$
\end{flushright}

\section{Lusternik-Schnirelmann ir-category}
Now we adapt the notion of Lusternik-Schnirelmann category that defined in [5], to the case of Lusternik-Schnirelmann ir-category.\\

\textbf{Definition 5.} The Lusternik-Schnirelmann ir-category of a space $X$ is the least integer $n$ such that there exists an open covering $U_{1}, \ldots,U_{n}$ of $X$ so that each $U_i$ is ir-contractible in $X$. We denote this by  $\stackrel{\rm ir}{\rm Cat}(X) = n$ and call such a covering $\{U_i\}$ ir-categorical. If no such integer exists, we write $\stackrel{\rm ir}{\rm Cat}(X) = \infty$. 

\vspace{1cm}

\textbf{Definition 6.} The space $X$ is said to be ir-path connected, if for each pair of points like $x,y$, there exist at least one ir-path from $x$ to $y$ or the reverse, from $y$ to $x$.  

\vspace{1cm}

\textbf{Notation.} For a ring $R$, we denote the set maximal ideals of $R$ by $Max(R)$.

\vspace{1cm}
\textbf{Theorem 8.} Let $R$ be a commutative ring. If the number of maximal ideals of $R$ is $k$ $\bigg ( \bigg |Max(R) \bigg | =k \bigg )$ then\\
\begin{center}
$\stackrel{\rm ir}{\rm Cat} \bigg (Spec(R) \bigg ) = k$
\end{center}

\vspace{1cm}
\textit{Proof.} consider $Max(R) = \bigg \{ M_1, \ldots , M_k \bigg \} $. As we know, the closed sets of $Spec(R)$ with Zariski topology are of the following form\\
\begin{center}
$ V(I) = \bigg \{ P \in Spec(R) \bigg | I \subseteq P \bigg \}
$
\end{center}
 Thus
\begin{center}
$ V(M_i) = \bigg \{ P \in Spec(R) \bigg | M_i \subseteq P \bigg \} = M_i
$
\end{center}
Now, we claim that  $\stackrel{\rm ir}{\rm Cat} \bigg (Spec(R) \bigg ) \ge k$. Suppose that  $U_{1}, \ldots,U_{m}$ is a ir-categorical cover for $Spec(R)$ and  $m < k$. If so, at least two maximal ideals like  $M_s$ and  $M_t$ placed in one open ir-contractible set of $Spec(R)$ like $U_r$. Since, $U_r$ is ir-contractible,  $\stackrel{\rm ir}{\rm Co} (U_r) = \bigcap_{P_{\alpha} \in U_r} \overline{\{P_{\alpha}\}} \neq \emptyset$. Also, as $M_t , M_s \in U_r$,  $\overline{\{M_t \}} = M_t$  and  $\overline{\{M_s \}} = M_s$, obviously $ \bigcap_{P_{\alpha} \in U_r} \overline{\{P_{\alpha}\}} = \emptyset$ and this is a contradiction. Thus $m \ge k$.\\
Now consider  $\hat{M_i} := (Max(R) \setminus M_i)$ and an ir-categorical cover for $Spec(R)$ in following
\begin{center}
$Spec(R) =  \bigcup_{i=1}^{k} \bigg(Spec(R) \setminus \hat{M_i} \bigg) $
\end{center}
Now we claim that for all $1 \le i \le k$ the open subset $W_i := Spec(R) \setminus \hat{M_i}$ of $Spec(R)$ is ir-contractible and  $\stackrel{\rm ir}{\rm Co} (W_i) = M_i$. The reason is the only maximal ideal in $W_i$ is $M_i$, and for each prime ideal in $W_i$ like  ${P_j}^{w_i}$ that $j \in J$, we have ${P_j}^{w_i} \subseteq M_i$, therefore $V(M_i) = M_i \subseteq V({P_j}^{w_i})$. Thus,  $M_i = \bigcap_{j \in J} \overline{ \bigg \{{P_j}^{w_i} \bigg \} }$ or equivalently  $M_i = \stackrel{\rm ir}{\rm Co} (W_i)$. Hence, all the $W_i$ for $1 \le i \le k$, are ir-contractible and give us an ir-categorical cover for $Spec(R)$.\\
Thus $\stackrel{\rm ir}{\rm Cat} \bigg (Spec(R) \bigg ) = k$.  

\begin{flushright}
$\square$
\end{flushright}

\vspace{1cm}

\textbf{Corollary 6.} If $R$ be a commutative and local ring, then $\stackrel{\rm ir}{\rm Cat} \bigg (Spec(R) \bigg ) = 1$, hence $Spec(R)$ is ir-contractible and $\stackrel{\rm ir}{\rm Co} \bigg (Spec(R) \bigg ) = M$, where $M$ is the only maximal ideal of $R$.
 
 \vspace{1cm}
 
 \textbf{Corollary 7.} Let $F$ be a field, then  $\stackrel{\rm ir}{\rm Co} (F) = \{0\}$, that means all the fields are ir-contractible.
 
 \vspace{1cm}
 
 \textbf{Theorem 9.} Let $X$ and $Y$ be topological spaces, then we have\\
 
\begin{center}
$\stackrel{\rm ir}{\rm Cat} (X \times Y) = \stackrel{\rm ir}{\rm Cat} (X) \ \times \stackrel{\rm ir}{\rm Cat}(Y)$
\end{center}

\vspace{1cm}

\textit{Proof.} Suppose that $W_1, \ldots , W_n$,  $V_1, \ldots , V_m$  and $(F_1 \times E_1), \ldots , (F_k \times E_k)$ be three ir-categorical covers respectively for  $X$, $Y$ and  $X \times Y$. Let $F_j$ be an arbitrary ir-contractible element of ir-categorical cover of $X$, therefore there exists a point like $x_j \in \ \stackrel{\rm ir}{\rm Co}(F_j)$. Also, there exists at least one  $1 \le i_0 \le n$ such that $x_j \in W_{i_0}$, thus  $x_j \in F_j \cap W_{i_0} $. Now we check all possible cases for this intersection:\\
1. If $F_j \cap W_{i_0} = W_{i_0} $, then  $x_j \in W_{i_0} \subset F_j$ and it contrasts with  $x_j \in  \ \stackrel{\rm ir}{\rm Co}(F_j)$. \\
2. If $F_j \cap W_{i_0} = \{O\} $  where  $O \subset F_j ,  W_{i_0}$, it contrasts with $x_j \in  \ \stackrel{\rm ir}{\rm Co}(F_j)$ similiar to the previous case.\\
3. If  $F_j \cap W_{i_0} = F_j $, we replace  $F_j$ with  $W_{i_0}$. \\
4. If $F_j \cap W_{i_0} = F_j = W_{i_0}$, we replace  $F_j$ with  $W_{i_0}$.\\

By repeating the same argument, we conclude that the elements of the ir-categorical cover of  $X \times Y$ are of the form  $W_i \times V_j$, where $1 \le i \le n$ and $1 \le j \le m$. Now we claim that $X \times Y = \bigcup_{j=1}^m \bigcup_{i=1}^n (W_i \times V_j)$ is an ir-categorical cover and denote it by $\mathcal{M}$. Assume that for a $1 \le i_0 \le n$ and a $1 \le j_0 \le m$,  $(W_{i_0} \times V_{j_0})$ is not in the $\mathcal{M}$. Since $(W_{i_0}$ and  $V_{j_0}$ are respectively in the ir-categorical covers of $X$ and $Y$, then they have some elements that are not in any other open ir-contractible subsets of their ir-categorical covers. Thus, $(W_{i_0} \times V_{j_0})$ has some elements that are not in any other open ir-contractible subsets of $\mathcal{M}$. This proves that $\mathcal{M}$ is an ir-categorical cover for $X \times Y$, and obviously $|\mathcal{M}| = n.m$ which implies the statement of the theorem.

\begin{flushright}
$\square$
\end{flushright}

\vspace{1cm}

\textbf{Theorem 10.} Let $X \subset  \ \stackrel{\rm ir}{\rm \mathbb{R}^n}$  with subspace topology of $\stackrel{\rm ir}{\rm \mathbb{R}^n}$, such that for all points $(a_1, \ldots , a_n) \in X$, there exists an unique point like  $(x_1, \ldots , x_n) \in X$ which for any $1 \le i \le n$, $a_i \le x_i$, then\\ \begin{center} $\stackrel{\rm ir}{\rm Co}(X) = (x_1 , \ldots , x_n)$ \end{center}

\vspace{1cm}

\textit{Proof.} We know that the closed sets of $X$ are in the following form\\ \begin{center} $\bigg (\prod_{i=1}^n [{\alpha}_i , \infty) \bigg ) \cap X = \prod_{i=1}^n [{\alpha}_i , x_i]$ \end{center} Thus it's obvious that $(x_1, \ldots , x_n)$  is in all closed sets of $X$, and therefore  $(x_1 , \ldots , x_n) \in \ \stackrel{\rm ir}{\rm Co}(X)$.\\
Now Suppose that $(x_1, \ldots , x_n) \neq (y_1, \ldots , y_n) \in \ \stackrel{\rm ir}{\rm Co}(X) $, thus $(y_1, \ldots , y_n)$ should be in all closed sets of $X$. But it is not in the closed set $\prod_{i=1}^n [{\frac{x_i + y_i}{2}}, x_i]$. Therefore $X$ is ir-contractible, only in $(x_1 , \ldots , x_n)$.

\begin{flushright}
$\square$
\end{flushright}

\vspace{1cm}

\textbf{Definition 7.} Let $X$ be a topological space and  $\gamma : \stackrel{\rm ir}{\rm I} \to X$ be an ir-path such that $\gamma(0) = x_0$ and $\gamma(1) = x_1$. The reverse ir-path of $\gamma$ is the ir-path $\overline{\gamma}:  \stackrel{\rm ir}{\rm I} \to X$ such that  $\overline{\gamma}(0) = x_1$ and  $\overline{\gamma}(1) = x_0$.

\vspace{1cm}

\textbf{Theorem 11.} If $X$ is a $T_0$ space and $\gamma : \stackrel{\rm ir}{\rm I} \to X$ is an ir-path, then $\gamma$ doesn't have a reverse ir-path.

\vspace{1cm}

\textit{Proof.} Suppose that $\gamma$ is an ir-path from $x$ to $y$ and $\overline{\gamma}$ is its reverse ir-path from $y$ to $x$. By applying the Theorem 2 we have   $y \in \overline{\{x\}}$ and $x \in \overline{\{y\}}$. Besides, since $X$ is $T_0$ we know that either there exists an open set   $U \subset X$  such that  $x \in U$ and  $y \in X \setminus U = V$, which since $V$ is a closed set containing $y$, then  $ \overline{\{y\}} \subseteq V$ and it means $x \in V$, which is contradiction, or there exists an open set  $U^{'} \subset X$ such that  $y \in U^{'}$ and $x \in X \setminus U^{'} =V^{'} $, which since $V^{'}$ is a closed set containing $x$, then $ \overline{\{x\}} \subseteq V^{'}$ and it means  $y \in V^{'}$, which is contradiction too. Thus $\gamma$ doesn't have any reverse ir-path.

\begin{flushright}
$\square$
\end{flushright}

\vspace{1cm}

\textbf{Theorem 12.} If $X$ be a $T_0$ and ir-contractible space, then  $\bigg |\stackrel{\rm ir}{\rm Co}(X) \bigg | = 1 $.

\vspace{1cm}

\textit{Proof.} Suppose that $x_1 , x_2 \in \stackrel{\rm ir}{\rm Co}(X)$, thus the only open set containing $x_1 \ , \ x_2$ is $X$. Also, since $X$ is $T_0$, then an open set like  $U \subset X$ should exists such that either $x_1 \in U$ and $x_2 \in X \setminus U$ or $x_2 \in U$ and  $x_1 \in X \setminus U$. Therefore $U$ must be equal to $X$ which is a contradiction. Thus, $\stackrel{\rm ir}{\rm Co}(X)$ is singleton.

\begin{flushright}
$\square$
\end{flushright}

\vspace{1cm}

\textbf{Proposition 3.} If  $X = W_1 \cup W_2 \cup \ldots \cup W_n$ is an ir-categorical cover for $X$ called $\mathcal{M}$, then for all $i,j \in \{1, \ldots , n\}$ ($i \neq j$), $\stackrel{\rm ir}{\rm Co}(W_i) \ \cap (W_j) = \emptyset$.

\vspace{1cm}

\textit{Proof.} Suppose that for one $i \ , \ j$ there exists a $x_0 \in X$ such that  $x_0  \in \ \stackrel{\rm ir}{\rm Co}(W_i) \ \cap W_j  $  and we have $W_i \ \cap W_j = U_0$. Since $\mathcal{M}$ is an ir-categorical cover, then $U_0 \neq W_i \ , \ W_j$. Thus  $U_0$  is an open set contained in $W_i$ and $W_j$ such that $x_0 \in U_0$. But it contrasts with  $x_0 \in \ \stackrel{\rm ir}{\rm Co}(W_i) $.

\begin{flushright}
$\square$
\end{flushright}

\vspace{1cm}

\textbf{Lemma 1.} Let $X = W_1 \cup W_2 \cup \ldots \cup W_n$ is an ir-categorical cover for $X$ called $\mathcal{W}$. If $\mathcal{V} = \{V_i \}_{i \in I}$ is an arbitrary open cover for $X$, then  $\mathcal{W}$ is a refinement of  $\mathcal{V}$.

\vspace{1cm}

\textit{Proof.} Let $x_i  \in \ \stackrel{\rm ir}{\rm Co}(W_i)$, therefore there exists a $V_{i_1} \ \in \mathcal{V}$ such that $x_i \ \in V_{i_1}$. Now if $W_i \ \cap V_{i_1} = U_i$ since $x_i \ \in U_i$, then $U_i = W_i$ and consequently  $W_i \subseteq V_{i_1}$. Thus $\mathcal{V} = \{V_i \}_{i \in I} = \{V_{i_j} \}_{j = 1}^n $ and  $\mathcal{W}$ is a refinement of  $\mathcal{V}$. 

\begin{flushright}
$\square$
\end{flushright}

\vspace{1cm}

\textbf{Corollary 8.} Let $X$ be a topological space, and $\mathcal{W}$ be its ir-categorical cover. Then $\mathcal{W}$ is the only refinement of $\mathcal{W}$.

\vspace{1cm}

\textbf{Lemma 2.} If $X$ is a topological space and $\stackrel{\rm ir}{\rm cat}(X) = n$, then $X$ has no cover with more than $n$ members.

\vspace{1cm}

\textit{Proof.} Suppose that $\mathcal{V} = \{V_j\}_{j=1}^{m}$ is an open cover for $X$ and  $X = W_1 \cup W_2 \cup \ldots \cup W_n$ be its ir-categorical cover. By using the proof of Lemma 1 we know that for every $1 \le i \le n$ there exists a  $V_{j_i} \in \mathcal{V}$ such that $W_i \subseteq V_{j_i}$. Therefore  $ X = \bigcup_{i=1}^n W_i \subseteq \bigcup_{j=1}^n V_{i_j}$ and hence  $\bigcup_{j=1}^n V_{i_j} = X$ which its cardinality is at most $n$.

\begin{flushright}
$\square$
\end{flushright}

\vspace{1cm}

\textbf{Theorem 13.}  If $X$ is a topological space and $\stackrel{\rm ir}{\rm cat}(X) = n$, then\\
\begin{center}
dim$(X) \ + 1 \le \ \stackrel{\rm ir}{\rm cat}(X) $
\end{center}

\vspace{1cm}

\textit{Proof.} By applying Lemma 2 we know that the order of all open covers is at most $m+1 = n$. Thus for the topological dimension of $X$ which is the minimum value of $m$, we have  dim$(X)+1 \le n$.

\begin{flushright}
$\square$
\end{flushright}

\vspace{1cm}

\textbf{Notation.} Let $X$ be a topological space. For $x,y \in X$, $x \preceq y$ denote that there exists an ir-path from $x$ to $y$.

\vspace{1cm}

\textbf{Proposition 4.} Let $X$ be a topological space. Then $\preceq$ is a quasiorder on $X$.

\vspace{1cm}

\textit{Proof.} Suppose that $\gamma_{x_0}$ is a constant ir-path such that for $t \in [0,1]$, $\gamma_{x_0}(t) = x_0$. Therefore there exists a constant ir-path from $x_0$ to $x_0$, and we can write $x_0 \preceq x_0$. Thus $\preceq$ is reflexive.\\
In order to prove that $\preceq$ is transitive, assume that $x \preceq y$ and $y \preceq z$. By using Theorem 2 we know that $y \in \overline{\{x\}}$ and $z \in \overline{\{y\}}$, and also $\overline{\{y\}} \subseteq \overline{\{x\}}$. Therefore $z \in \overline{\{x\}} $ and it means $x \preceq z$.
\begin{flushright}
$\square$
\end{flushright}

\vspace{1cm}

\textbf{Corollary 9.} If $X$ be a $T_0$ space, then $\preceq$ is a partial order over $X$.

\vspace{1cm}

\textit{Proof.} Follows from Theorem 11.

\begin{flushright}
$\square$
\end{flushright}

\section{ir-homotopy equivalence}

\textbf{Definition 8.} Let $X$ and $Y$ be two topological spaces. A map $f: X \to Y$ is called an ir-homotopy equivalence if there exists a map $g: Y \to X$ such that $1_X \stackrel{\rm ir}{\rm \cong} f \circ g$ and $1_Y \stackrel{\rm ir}{\rm \cong} g \circ f$. If an ir-homotopy equivalence exists, the spaces $X$ and $Y$ are said to be ir-homotopy equivalent.

\vspace{1cm}

\textbf{Theorem 14.} Let $X$ and $Y$ be ir-homotopy equivalent. If $X$ is ir-contractible then $Y$ is ir-contractible.

\vspace{1cm}

\textit{Proof.} We know that there exist maps $f: X \to Y$ and $g: Y \to X$ such that $1_X \stackrel{\rm ir}{\rm \cong} f \circ g$ and $1_Y \stackrel{\rm ir}{\rm \cong} g \circ f$. Suppose that $x_0 \in \ \stackrel{\rm ir}{\rm Co} (X)$, we show that $f(x_0) \in \ \stackrel{\rm ir}{\rm Co}(Y)$. For proving we should show that, if an open subset $V$ of $Y$ contains $f(x_0)$, then $V = Y$. Since $f(x_0) \in V$ then $x_0 \in f^{-1}(V)$, and also because $f^{-1}(V)$ is open, therefore $f^{-1}(V) = X$ and consequently $f(X) \subseteq V$. Thus for all arbitrary $y \in Y$, it's clear that $f(g(y)) \in V$. \\
Now consider the ir-homotopy $G: Y \times  \stackrel{\rm ir}{\rm I} \to Y$ where for all $y \in Y$, $G(y,0) = y$ and $G(y,1) = f(g(y))$. For an arbitrary $y$, we define an ir-path $\gamma : \stackrel{\rm ir}{\rm I} \to Y$ with the formula $\gamma(t) = G(y,t)$. Since $\gamma(0) = y$ and $\gamma(1) = f(g(y))$ by using Theorem 2, $f(g(y)) \in \overline{\{y\}}$. Thus the open neighborhood $V$ of $f(g(y))$ contains $y$.

\begin{flushright}
$\square$
\end{flushright}

\vspace{1cm}

\textbf{Theorem 15.} Let $X$ and $Y$ be ir-homotopy equivalent. Then  $\stackrel{\rm ir}{\rm cat (X)} = \stackrel{\rm ir}{\rm cat (Y)}$.

\vspace{1cm}

\textit{Proof.} Since $X$ and $Y$ are ir-homotopy equivalent, we have continuous maps $f: X \to Y$  and $g: Y \to X$ such that\\

\begin{center}
$
1_X \stackrel{\rm ir}{\rm \cong} g \circ f 
$
\end{center}

\begin{center}
$
1_Y \stackrel{\rm ir}{\rm \cong} f \circ g
$
\end{center}

Assume that  $\stackrel{\rm ir}{\rm cat (X)} = n$, and  $X = w_1 \cup w_2 \cup \ldots \cup w_n$ is an ir-categorical cover for $X$ and  $x_i \in \stackrel{\rm ir}{\rm Co}(w_i)$ we show that $f(x_i) \in \stackrel{\rm ir}{\rm Co}(g^{-1}(w_i))$. Now since  $1_X \stackrel{\rm ir}{\rm \cong} g \circ f$, we have the following ir-homotopy\\

\begin{center}
$
F: X \times  \stackrel{\rm ir}{\rm I} \to X$ \\
$F(x,0) = x$ \\
$F(x,1) = g( f(x))
$
\end{center}

By putting $x := x_i$ we have an ir-path as the following\\
\begin{center}
$F_{x_i} = \gamma :  \stackrel{\rm ir}{\rm I} \to X $\\
$\gamma (0) = x_i$ \\
$\gamma (1) = g ( f (x_i))$
\end{center}

Since $g( f (x_i)) \in X$ thus there exists a  $j \in \{1, \ldots, n\}$ such that $g( f (x_i)) \in w_j$, now if  $w_j \neq w_i$ since $\gamma^{-1}(w_j) = \stackrel{\rm ir}{\rm I} $, therefore  $x_i \in w_j $ that by using Proposition 29 it's a contradiction. Thus  $g( f (x_i)) \in w_i$ or equivalently $f(x_i) \in g^{-1}(w_i)$. It is sufficient to show that if arbitrary open subset $V$ of $g^{-1}(w_i)$ contains  $f(x_i)$, then  $V = g^{-1}(w_i)$. Now from  $f(x_i) \in V$ we have $x_i \in f^{-1}(V)$ and since  $x_i \in \stackrel{\rm ir}{\rm Co}(w_i)$ then $w_i \subseteq  f^{-1}(V)$ and therefore $f(w_i) \subseteq V$. Also, we know $g(g^{-1}(w_i)) \subseteq w_i$, consequently  $f(g(g^{-1}(w_i))) \subseteq f(w_i) \subseteq V$, thus for each $y \in g^{-1}(w_i)$, $f(g(y)) \in V$. On the other hand, we have the following ir-homotopy\\

\begin{center}
$
G: Y \times  \stackrel{\rm ir}{\rm I} \to Y$ \\
$G(y,0) = y$ \\
$G(y,1) = f(g(y))
$
\end{center}

For arbitrary $y$ we can define the following ir-path\\

\begin{center}
$G_y = \sigma :  \stackrel{\rm ir}{\rm I} \to Y $\\
$\sigma (0) = y$ \\
$\sigma (1) = f(g(y))$
\end{center}

Now for each $y \in g^{-1}(w_i)$, the open neighborhood $V$ contains  $f(g(y))$ and therefore  $y \in V$ which is clear from  $f(g(y)) \in \overline{\{y\}}$. Thus   $V = g^{-1}(w_i)$ which gives $f(x_i) \in \stackrel{\rm ir}{\rm Co}(g^{-1}(w_i))$.\\
Now since  $Y = g^{-1} (X) = \bigcup_{i=1}^{n} g^{-1} (w_i) $ thus  $\bigcup_{i=1}^{n} g^{-1} (w_i)$ is an open cover for $Y$ such that its elements are ir-contractible, hence $\stackrel{\rm ir}{\rm cat} (Y) \ \le \ \stackrel{\rm ir}{\rm cat} (X)$.\\

By repeating the same argument, we conclude that $\stackrel{\rm ir}{\rm cat} (X) \ \le \ \stackrel{\rm ir}{\rm cat} (Y)$. Thus we proved the statement, and $\stackrel{\rm ir}{\rm cat} (X) = \stackrel{\rm ir}{\rm cat} (Y)$.

\begin{flushright}
$\square$
\end{flushright}

\vspace{1cm}

\textbf{Remark.} Theorem 14 is a corollary of the Theorem 15.

\section*{} \label{bibsection}

\bibliographystyle{plain}
\bibliography{template}

\end{document}